\documentclass[10pt]{amsart}
\usepackage{amssymb, amsthm, amsmath}
\usepackage{latexsym}
\usepackage{epsfig}

\theoremstyle{plain}

\newtheorem{thm}{Theorem}

\theoremstyle{remark}

\theoremstyle{definition}

\def\A{{\mathcal A}}

\def\Z{{\mathbb Z}}

\def\R{{\mathbb R}}
\def\C{{\mathbb C}}

\def\l{{\lambda}}
\def\m{{\mu}}
\def\n{{\nu}}

\def\sh{{\mathfrak s}}

\def\sR{{\mbox{R}}}
\def\cR{{\mathcal R}}
\def\a{{\alpha}}
\def\b{{\beta}}

\def\s{{\sigma}}
\def\t{{\tau}}
\def\om{{\omega}}
\def\th{{\theta}}
\def\Om{{\Omega}}

\newcommand{\dis}{\displaystyle}

\newcommand{\gequ}{\geqslant}
\newcommand{\lequ}{\leqslant}
\newcommand{\lra}{\longrightarrow}
\newcommand{\hra}{\hookrightarrow}
\newcommand{\ra}{\rightarrow}

\newcommand{\dbar}{\overline{\partial}}

\newcommand{\Pol}{\mbox{Pol}}
\newcommand{\End}{\mbox{End}}

\newcommand{\Ker}{\mbox{Ker}}
\newcommand{\Tr}{\mbox{Tr}}
\newcommand{\id}{\mbox{id}}
\newcommand{\ov}{\overline}
\newcommand{\noin}{\noindent}
\newcommand{\wt}{\widetilde}

\newcommand{\med}{\medskip}

\newcommand{\lie}[1]{{\mathfrak #1}}

\newcommand{\gl}{\lie{g}\lie{l}}

\begin{document}

\title[From representation theory to Schubert calculus]
{The connection between representation theory and Schubert calculus}
\author{Harry Tamvakis}
\date{December 20, 2004}
\thanks{The author was supported in part by National Science 
Foundation Grant DMS-0296023}
\address{Department of Mathematics, Brandeis University - MS 050,
P. O. Box 9110, Waltham, MA
02454-9110, USA}
\email{harryt@brandeis.edu}

\maketitle

\section{Introduction}
\label{intro}

\noin

The irreducible polynomial representations of the general linear group
$GL_n(\C)$ are parametrized by integer partitions $\l$ with at most
$n$ parts. Given any two such representations $V^{\l}$ and $V^{\m}$,
one has a decomposition of the tensor product
\begin{equation}
\label{gltensor}
V^{\l}\otimes V^{\m} = \sum_{\n} c_{\l\m}^{\n} \,V^{\n}
\end{equation}
into irreducible representations $V^{\nu}$ of $GL_n$. 

Let $G(m,n)$ denote the Grassmannian of complex $m$-dimensional linear
subspaces of $\C^{m+n}$. The cohomology ring of $G(m,n)$ has a natural 
geometric basis of Schubert classes $\s_{\l}$, 
%again parametrized by partitions. 
and there is a cup product decomposition 
\begin{equation}
\label{gmncup}
\s_{\l}\cdot\s_{\m} = \sum_{\n} c_{\l\m}^{\n}\, \s_{\n}.
\end{equation}
The structure constants $c_{\l\m}^{\n}$ determine the classical 
Schubert calculus on $G(m,n)$.

It has been known for some time that the integers $c_{\l\m}^{\n}$ in
formulas (\ref{gltensor}) and (\ref{gmncup}) coincide, as long as the
meaningless Schubert classes in (\ref{gmncup}) are interpreted as
zero. Following the work of Giambelli \cite{G} \cite{G2}, this is
proved formally by relating both products to the multiplication of
Schur $S$-polynomials; a precise argument along these lines was given
by Lesieur \cite{Les}.
%Lesieur \cite{Les}, this is 
%proved formally by relating both products
%to the multiplication of Schur $S$-polynomials
%(see e.g.\ \cite[\S 9.4]{F1}).
It is natural to ask for a more direct, conceptual explanation of this 
fact. This question has appeared 
%occasionally 
every so often in print; some recent examples are
\cite[\S 6.2]{F2} and \cite[\S 1]{Len}.

Our aim here is to describe a direct and natural connection between
the representation theory of $GL_n$ and the Schubert calculus, which
goes via the Chern-Weil theory of characteristic classes. Indeed,
since the Grassmannian is a universal carrier for the Chern classes of
principal $GL_n$-bundles, it is not so surprising that the cohomology
ring of $G(m,n)$ is related to the representation ring of $GL_n$.  
From this point of view, we can also understand why a
result of this sort fails to hold for other types of Lie groups: what
makes $GL_n$ special is the fact that it sits naturally as a dense
open subset of its own Lie algebra (see Sec.\ \ref{rsf}).
The key observation is that the Chern-Weil homomorphism 
extends to a ring homomorphism from the (polynomial) 
representation ring $\sR_+(GL_n)$ to $H^*(G(m,n))$, which 
sends the natural basis elements of the first ring to the Schubert classes.

The relation between Schubert calculus and the multiplication of Schur
polynomials has been investigated before by Horrocks \cite{H} 
and Carrell \cite{Ca}. Although the approach in \cite{H}  is closest
to the one here, the main ideas go back to the fundamental
works of Chern \cite{Ch1}, Weil \cite{W}, and H.\ Cartan 
\cite{Car}.  We provide
an exposition where the various ingredients from representation theory,
differential geometry, topology of fiber bundles, and Schubert calculus
are each presented in turn. In Sec.\ \ref{prob}, we apply Grothendieck's
construction of the Chern classes of Lie group representations to look
for an analogue of these results in the other Lie types.

%I thank Bill Fulton for asking the right question and encouraging me to 
%write this article.
%I am also grateful to Frank Sottile for similar encouragement
%and for bringing Horrocks' paper \cite{H} to my attention, Alain
%Lascoux for information on Giambelli's work, 
%and Danny Ruberman for useful discussions.
%Finally, a recent preprint of Beauville \cite{Be} helped clarify the relation 
%between the representation rings of  Lie groups and the 
%characteristic homomorphism.

I thank Bill Fulton and Frank Sottile for encouraging me to write 
this article,  and also Arnaud Beauville, Allen Knutson, 
Alain Lascoux, and Danny Ruberman
for useful discussions and email exchanges.

\section{Representations and Schur functors}
\label{rsf}

We are concerned here with the polynomial representations of the
general linear group $GL_n(\C)$.  A matrix representation 
$\pi:GL_n\ra GL_N$ of $GL_n$ is 
{\em polynomial} if the entries of  $\pi(g)$ are polynomials in the entries
of $g\in GL_n$. The character of $\pi$ is the function $\chi:GL_n \ra 
\C$ defined by $\chi(g)=\Tr(\pi(g))$. The polynomial representation
ring $\sR=\sR_+(GL_n)$ is the
$\R$-algebra generated by the polynomial characters of $GL_n$. The
ring $\sR$ may be identified with the real vector space spanned by the 
irreducible polynomial
$GL_n$-representations (up to isomorphism), with the ring structure 
given by the tensor product. We use real instead of integer coefficients here
because the Chern-Weil construction in the sequel will involve de Rham
cohomology groups.

The group $GL_n$ acts by conjugation on its Lie algebra $\gl_n$ (the
space of all $n\times n$ matrices). This induces a $GL_n$ action on
the ring $\Pol(\gl_n)$ of polynomial functions on $\gl_n$ with real
coefficients; we denote by $\Pol(\gl_n)^{GL_n}$ the corresponding ring
of invariants. Since $GL_n$ is a dense open subset of $\gl_n$, any
character $\chi\in \sR$ extends to a unique element of
$\Pol(\gl_n)^{GL_n}$. Conversely, for any invariant polynomial
$f\in\Pol(\lie{g}\lie{l}_n)^{GL_n}$, the restriction of $f$ to $GL_n$
is clearly a polynomial class function, and hence an element of
$\sR$. Thus we obtain a canonical ring isomorphism
\begin{equation}
\label{isom1}
\phi: \sR_+(GL_n)\ra \Pol(\gl_n)^{GL_n}.
\end{equation}
In contrast, there is no satisfactory analogue of the morphism $\phi$
for the other types of Lie groups (see Sec.\ \ref{prob}). Note that
there is already the problem of defining `polynomial representations'
for a general complex Lie group.\footnote{For a connected reductive
complex Lie group $G$ with Lie algebra $\lie{g}$, Knutson
suggests to define a ring of polynomial representations of $G$ as
the image of the injective map $\Pol(\lie{g})^G \ra \Pol(G)^G$ which
is induced by pullback along a generalized Cayley transform
$G\ra \lie{g}$ (see \cite{KM}).}
 
Following Schur \cite{S1} \cite{S2}, the irreducible polynomial
representations of $GL_n$ are parametrized by integer partitions
$\l=(\l_1\gequ \l_2\gequ\cdots \gequ \l_n\gequ 0)$ of length
(i.e. number of nonzero parts $\l_i$) at most $n$. For each partition
$\l$ there is a {\em Schur functor} $\sh_{\l} : {\mathfrak V} \ra
{\mathfrak V}$, where ${\mathfrak V}$ denotes the category of finite
dimensional $\C$-vector spaces and $\C$-linear maps. If $V=\C^n$ is
the standard representation of $GL_n$, then the irreducible
representation corresponding to $\l$ is the Schur module
$V^{\l}=\sh_{\l}(V)$. In the language of Lie theory, $V^{\l}$ is the
$GL_n$-representation with highest weight $\l$. The character of
$V^{\l}$ is a {\em Schur polynomial} in the eigenvalues of $g\in
GL_n$.

For completeness, we briefly describe Weyl's construction of $\sh_{\l}(V)$,
for any complex vector space $V$. 
First, we identify the partition $\l$ with its
Young diagram of boxes; this is an array of $p=\sum \l_i$ boxes arranged 
in left-justified rows, with $\l_i$ boxes in the $i$th row. Number the 
boxes in $\l$ with the integers $1,\ldots ,p$ in order going from left
to right and  top to bottom. The resulting standard tableaux $T$
is illustrated below on the Young diagram of $\l=(4,3,2)$.

\medskip

$$\includegraphics[scale=.32]{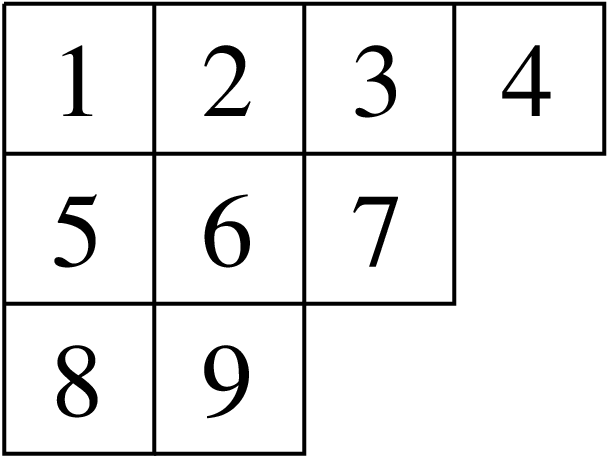}$$

\medskip
\noin
Let $R$ (respectively $C$) denote the subgroup of the symmetric group 
$S_p$ consisting of elements which permute the entries of each row 
(respectively column) of $T$ among themselves. Consider the elements
\[
a_{\l}=\sum_{u\in C} \mathrm{sgn}(u)\, u , \quad\qquad b_{\l}=\sum_{v\in R}v
\]
of the group algebra $\C[S_p]$, and define the {\em Young symmetrizer}
$c_{\l}=a_{\l}b_{\l}$. 

For any vector space $V$, $S_p$ acts on the right 
on the $p$-fold 
tensor product $M=V^{\otimes p}$ by permuting the factors, and this action 
commutes with the left action of $GL(V)$ on $M$. The {\em Schur module}
$\sh_{\l}(V)$ is defined as the image of the map $M\ra M$ that is right
multiplication by $c_{\l}$. This construction is functorial, in the sense that
a linear map $f:V\ra W$ of vector spaces determines a linear map
$\sh_{\l}(f):\sh_{\l}(V)\ra\sh_{\l}(W)$, with $\sh_{\l}(f_1\circ f_2)=
\sh_{\l}(f_1)\circ \sh_{\l}(f_2)$ and $\sh_{\l}(\mathrm{id}_V)=
\mathrm{id}_{\sh_{\l}(V)}$. 
In passing, we note that the {\em Specht modules}
$\C[S_p] c_{\l}$ form a complete set of 
irreducible representations of the symmetric group $S_p$, 
as $\l$ varies over all partitions of $p$. For more 
information on Schur modules the reader may consult e.g.\
\cite{FH}, \cite[\S 8]{F1}, \cite{Gr}, and \cite{Ma}. 
%and \cite[\S I, Appendix A]{M}.

\medskip
\noin
{\bf Example 1.} Let $p=2$ and write $S_2=\{1,\sigma\}$. The two
relevant partitions $\l$ are $(2)$ and $(1,1)$, with respective
Young symmetrizers  $c_2=1+\sigma$ and 
$c_{1,1}=1-\sigma$.  For any vector space $V$, there is a decomposition
$V\otimes V = \mathrm{Sym}^2V\oplus \wedge^2V$, and we see that 
$\sh_2(V)=\mathrm{Sym}^2V$ and $\sh_{1,1}(V)=\wedge^2V$.
\medskip

Given an $n\times n$ matrix $A=\{a_{ij}\}$ of indeterminates, we
define the {\em Schur matrix} to be $\sh_{\l}(A)$. This is a square matrix 
of size $\dim(V^{\l})$ whose entries are polynomials in the $a_{ij}$
with integer coefficients.  Moreover, the map $\phi$ in (\ref{isom1})
sends $V^{\l}$ to the invariant polynomial $A\mapsto
\Tr(\sh_{\l}(A))$.  It is surprising that
%more than $100$ years after Schur's thesis, there is 
no explicit formula for the entries of $\sh_{\l}(A)$ is known, except
in special cases. In general, there are several algorithms available
for this computation, some rather classical (see \cite{CLL},
\cite{Cl}, \cite{DKR}, and \cite{GK} for a sample).

\medskip
\noin
{\bf Example 2.} 
The fundamental representations of $GL_n$ are the exterior powers 
$\wedge^kV$, for $V=\C^n$ and $k=1,\ldots,n$, which correspond to 
the partitions $\l=(1^k)$. The Schur matrix $\wedge^kA$ has order 
${n \choose k}$, and its rows (resp.\ columns) are indexed by the $k$-element 
subsets of the $n$ rows (resp.\ columns) of $A$. The entries of $\wedge^kA$ 
are the determinants of the $k\times k$ minors in $A$. The corresponding 
invariant polynomials for $k=1$ and $k=n$ are given by
$\Tr(A)$ and $\det(A)$, respectively. 

\medskip
 
If $A$ is a diagonal matrix with eigenvalues $x_1,\ldots,x_n$, then
$\Tr(\sh_{\l}(A))$ is the Schur polynomial $s_{\l}(x_1,\ldots,x_n)$.
There is much to say about these important polynomials (see e.g.\
\cite[\S I]{M}), but we shall refrain from doing so because they are
not used in the sequel. We simply note here that since the Schur
polynomials form a basis of the ring of symmetric polynomials in
$x_1,\ldots,x_n$, it follows that the rings in (\ref{isom1}) are both
isomorphic to $\R[x_1,\ldots,x_n]^{S_n}$.

\section{The Chern-Weil homomorphism}
\label{cw}

The homomorphism which is the subject of this section is a fundamental tool 
in the theory of characteristic classes. Our main references are the monographs
\cite{Du}, \cite{GHV}, and \cite{GS}, all of which contain detailed
expositions of Chern-Weil theory and the related work of H.\ Cartan.

Consider a principal $GL_n$-bundle $P\ra M$ over a differentiable
manifold $M$, and choose an open cover $\{U_{\a}\}$ of $M$ which
trivializes $P$. The transition functions for $P$ with respect to this
cover are a system of morphisms $g_{\a\b}:U_{\a}\cap U_{\b}\ra GL_n$
satisfying the cocycle condition $g_{\a\b}g_{\b\gamma}g_{\gamma\a}=1$
on the intersections $U_{\a}\cap U_{\b}\cap U_{\gamma}$.  We may
identify $P\ra M$ with the rank $n$ complex vector bundle with the
same transition functions $g_{\a\b}$; in other words, with the
associated vector bundle for the standard representation of $GL_n$ on
$\C^n$.

For each $g\in GL_n$, let $R_g:P\ra P$ denote the right action of
$GL_n$ on $P$, and for each $x\in P$, let $v_x:\gl_n\ra T_xP$ be the
differential of the map $g\mapsto xg$.  A {\em connection} on $P$ is a
smooth $\gl_n$-valued $1$-form $\th\in {\mathcal A}^1(P,\gl_n)$
satisfying (i) $\th_x v_x = \id$ for all $x\in P$, and (ii)
$R_g^*\th=\mbox{Ad}(g^{-1})\th$, for all $g\in GL_n$. Note that (ii)
asserts that $\th$ is equivariant with respect to the adjoint
representation of $GL_n$. Connections always exist, provided the
base manifold $M$ admits partitions of unity.

The {\em curvature} $\Theta\in {\mathcal A}^2(P,\gl_n)$ of the
connection $\th$ is defined by the Maurer-Cartan
%\footnote{This is \'Elie Cartan, the father of Henri Cartan.} 
structure equation $\Theta=d\th+\frac{1}{2}[\th,\th]$ (this is \'Elie
Cartan, the father of Henri Cartan).  A basic theorem of Weil states
that for any homogeneous polynomial $f\in \Pol(\gl_n)^{GL_n}$ of
degree $k$, (i) the differential form $f(\frac{i}{2\pi}\Theta)$ is
closed in ${\mathcal A}^{2k}(M,\R)$, and (ii) its de Rham cohomology
class does not depend on the connection $\th$.  The resulting map
\[
\Pol(\gl_n)^{GL_n} \ra H^*(M,\R); \qquad
f \mapsto f(\frac{i}{2\pi}\Theta)
\]
is an algebra homomorphism called the {\em Chern-Weil
homomorphism}. We will continue to use cohomology with real
coefficients unless otherwise indicated.

In topology, one constructs a contractible space $EGL_n$ on which
$GL_n$ acts freely, and with quotient equal to the classifying space
$BGL_n$.  Every principal $GL_n$-bundle over $M$ is a pullback of the
universal bundle $EGL_n \ra BGL_n$; in this way, the isomorphism
classes of principal $GL_n$-bundles over $M$ are in one-to-one
correspondence with the homotopy classes of maps from $M$ to $BGL_n$.
Moreover, the previous definitions of connection and curvature have
universal analogues, and there is a homomorphism
$\Pol(\gl_n)^{GL_n}\ra H^*(BGL_n)$ (see \cite{NR} or \cite[\S 5,6]{Du}
for a detailed construction).  Since the quotient $GL_n/U(n)$ of
$GL_n$ by the unitary group $U(n)$ is diffeomorphic to a Euclidean
space, the inclusion $U(n)\hra GL_n$ induces an isomorphism
$H^*(BGL_n) \ra H^*(BU(n))$. We deduce that there is a universal
Chern-Weil homomorphism
\begin{equation}
\label{ucw}
\psi:\Pol(\gl_n)^{GL_n}\ra H^*(BU(n)).
\end{equation}
H. Cartan \cite{Car} proved that the map $\psi$ is an isomorphism of 
polynomial rings.
The cohomology of $BU(n)$ is thus identified with the ring of characteristic
classes of principal $GL_n$-bundles (or complex vector bundles). If the 
polynomial $f\in \Pol(\gl_n)^{GL_n}$ has integer coefficients, then its
image under the Chern-Weil map (\ref{ucw}) lies in $H^*(BU(n),\Z)$
(which we identify with its image in $H^*(BU(n),\R)$).

\med \noin {\bf Example 3.} Let $E\ra M$ be a rank $n$ complex vector
bundle over $M$.  The Chern-Weil homomorphism sends the invariant
polynomial $A\mapsto \Tr(\wedge^kA)$ of the previous section to the
$k$-th Chern class $c_k(E)\in H^{2k}(M,\Z)$.

\section{Schubert varieties and Schubert forms}
\label{svcc}

The cohomology ring of the Grassmannian $X=G(m,n)$ has a basis of 
Schubert classes $\s_{\l}$, one for each partition $\l$ whose Young diagram
is contained in an $n\times m$ rectangle (equivalently, $\l$ is a partition of
length at most $n$ with $\l_1\lequ m$). 
The Schubert class $\s_{\l}$ may be defined 
using the Poincar\'e duality isomorphism between cohomology and 
homology, as follows: $\s_{\l}$ is 
the element of $H^*(X)$  whose cap product with the fundamental 
class of $X$ is the homology class of a 
{\em Schubert variety} $X_{\l}$, 
described below. If the diagram of $\l$ does not
fit in the above rectangle, then we set $\s_{\l}=0$.

To define $X_{\l}$, consider the fixed complete flag of subspaces
\[
0=F_0\subset F_1\subset \cdots \subset F_N=\C^N
\]
where $N=m+n$ and $F_j=\C^j\times\{0\}\subset \C^N$, for $0\lequ j \lequ N$.
Let $\l'$ denote the conjugate 
partition, whose diagram is the transpose of the diagram of $\l$. 
The variety $X_{\l}$ consists of those $m$-dimensional subspaces 
$U$ such that $\dim(U\cap F_{n+i-\l_i'})\gequ i$, for $1\lequ i\lequ m$.
The complex codimension of $X_{\l}$ in $X$ is
given by the weight $|\l|=\sum \l_i$ of $\l$.
Note that our indexing convention for Schubert varieties is the transpose 
of the usual one, as found for example in \cite[\S 4]{F2}.

There is a tautological short exact sequence of vector bundles
\begin{equation}
\label{uses}
0\ra S \ra E \ra Q \ra 0
\end{equation}
over $X$, where $E=X\times \C^N$ is the trivial bundle of rank $N$ over $X$,
$S$ denotes the tautological
rank $m$ subbundle of $E$, and $Q$ the quotient bundle. The total space of 
$S$ is the submanifold of the product $X\times\C^N$ consisting of
those pairs $([U],v)$ with $v\in U$, and the projection $X\times \C^N\ra X$ 
restricts to the map from $S$ to $X$. Observe that the fiber of $S$
(resp.\ $Q$) over a point $[U]$ in $X$ which corresponds to the subspace 
$U\subset \C^N$ is given by $U$ itself (resp.\ $\C^N/U)$. From now on
it will be convenient to work with the quotient bundle $Q$ and to think
of $X=G(m,n)$ as parametrizing rank $n$ quotients of $\C^N$. 

According to \cite[\S 23]{BT} (see also \cite[\S 14]{MS} and \cite[\S
18]{Hu}), the {\em infinite Grassmannian}
$G(\infty,n)=\lim_{m\ra\infty}G(m,n)$ provides a model for the
classifying space $BU(n)$ of the unitary group. There is a 
natural sequence of inclusions
\[
\cdots \hra G(m-1,n) \hra G(m,n) \hra G(m+1,n)\hra \cdots
\]
and one defines the space 
$G(\infty,n)$ as the union of all the finite Grassmannians
$G(m,n)$ over $m\gequ 1$, with the inductive topology.
%The space
%$G(\infty,n)$ is defined by a telescoping construction applied to the
%natural sequence of inclusions
%\[
%\cdots \hra G(m-1,n) \hra G(m,n) \hra G(m+1,n)\hra \cdots .
%\]
The inclusion $G(m,n)\hra G(\infty,n)=BU(n)$ induces a surjection
\[
\zeta :H^*(BU(n))\ra H^*(G(m,n)).
\] 
Moreover, there is a universal rank $n$ quotient bundle $Q\ra
G(\infty,n)$ which corresponds to the universal bundle $EU(n)\ra
BU(n)$.

We let $\rho=\zeta\circ\psi\circ \phi$ denote the
composite of the three maps
\begin{equation}
\label{rho}
\sR_+(GL_n)\stackrel{\phi}\lra
\Pol(\gl_n)^{GL_n} \stackrel{\psi}\lra
H^*(BU(n))\stackrel{\zeta}\lra H^*(G(m,n)).
\end{equation}
The connection between the representation ring of $GL_n$ and the 
cohomology of $G(m,n)$ is exhibited in the following 
result.\footnote{The pun here and in the title of this paper is intended.}

\begin{thm}
\label{mainthm}
For every $\l$, the ring homomorphism $\rho \colon ${\em R}$_+(GL_n)
\ra H^*(G(m,n))$ maps the class of the irreducible representation
$V^{\l}$ to the Schubert class $\s_{\l}$.
\end{thm}
 
\med

As a prelude to the proof of Theorem \ref{mainthm}, we will describe
an equivalent method of constructing the morphism $\rho$. This
involves putting a $U(N)$-invariant hermitian metric on the vector
bundle $Q\ra X$ and using the canonical induced linear connection to
define {\em Schubert forms} $\Om_{\l}$ which represent the Schubert
classes in the de Rham cohomology of $X$. We will denote by $\A^k(X)$
(respectively $\A^k(X,Q)$) the real vector space of $C^{\infty}$
$k$-forms on $X$ (respectively, $Q$-valued $k$-forms on $X$).  A
connection on $Q$ is a $\C$-linear map $D: \A^0(X,Q)\ra \A^1(X,Q)$
such that
\begin{equation}
\label{coneq}
D(f\cdot s) = df\otimes s + f\cdot ds
\end{equation}
for all functions $f\in \A^0(X)$ and sections 
$s\in \A^0(X, Q)$. 
The type decomposition 
$\A^1(X,Q)=\A^{1,0}(X,Q)\bigoplus \A^{0,1}(X,Q)$ 
of differential forms   induces a decomposition 
$D=D^{1,0}+D^{0,1}$ of each connection $D$ on $Q$.

The standard hermitian metric on $\C^N$  gives a metric on the
trivial vector bundle in (\ref{uses}) and induces a metric on the subbundle
$S$. One obtains a hermitian metric $h$ on the quotient bundle $Q$ by
identifying  it with the orthogonal complement of $S$.  
The metric $h$ induces a unique connection $D=D(Q,h)$
such that $D^{0,1}=\overline{\partial}_Q$ and $D$ is {\em unitary}, i.e. 
\[
\dis
d\,h(s,t)=h(Ds,t)+h(s,Dt), 
\mbox{ for all } s,t\in \A^0(X,Q).
\]
The connection $D$ is called the {\em hermitian
holomorphic connection} of $(Q,h)$. 
We extend $D$ to $Q$-valued  forms by using the Leibnitz rule
\[
D(\om\otimes s)=d\om \otimes s + (-1)^{\deg \om}\om\otimes Ds.
\]

The curvature of $D$ is the composite
\[
\wt{\Om}=D^2:\A^0(X,Q)\ra \A^2(X,Q).
\]
By applying (\ref{coneq}) twice we compute that $\wt{\Om}(f\cdot s)=
f\cdot \wt{\Om}(s)$, hence
the map $\wt{\Om}$ is $\A^0(X)$-linear. We deduce that
$\wt{\Om}\in\A^2(X,\End(Q))$. In fact, 
$\wt{\Om}=D^{1,1}\in
\A^{1,1}(X,\End(Q))$, because $D^{0,2}=\dbar^2_Q=0$, so $D^{2,0}$ also
vanishes by unitarity. It follows that locally, 
we can identify $\wt{\Om}$ with an $n\times n$ matrix of $(1,1)$-forms 
on $X$.

Let $\Om=\frac{i}{2\pi}\wt{\Om}$. For each partition $\l$, define the Schubert
form $\Om_{\l}=\Tr(\sh_{\l}(\Om))$; this is a closed form of type $(|\l|,|\l|)$
on $X$. Since the hermitian vector bundle $(Q,h)$ is $U(N)$-equivariant
(for the natural $U(N)$ action on $X$),
the Schubert forms $\Om_{\l}$ are $U(N)$-invariant. The class
of  $\Om_{\l}$ in the de Rham cohomology group $H^{2|\l|}(X)$
coincides with the image $\rho(V^{\l})$ in Theorem \ref{mainthm}. One
has an equation
\[
\Tr(\sh_{\l}(\Om))\wedge
\Tr(\sh_{\mu}(\Om)) = 
\Tr(\sh_{\l}(\Om)\otimes
\sh_{\mu}(\Om)) = 
\sum_{\n} c_{\l\m}^{\n} \Tr(\sh_{\n}(\Om))
\]
of differential forms on $X$, and hence a formula
\[
[\Om_{\l}]\cdot [\Om_{\m}]=[\Om_{\l}\wedge \Om_{\m}]
=\sum_{\n} c_{\l\m}^{\n} [\Om_{\n}]
\]
in $H^*(X)$, with the constants $c_{\l\m}^{\n}$ defined as in (\ref{gltensor}). 

To prove the theorem, we must show that $[\Om_{\l}]$ is equal to the
Schubert class $\s_{\l}$. The link between Schubert classes and
$GL_n$-modules is then evident: each Schubert class $\s_{\l}$
is represented  in de Rham cohomology by a unique\footnote{Since 
the Grassmannian $X$ is a hermitian symmetric space,
the $U(N)$-invariant forms are harmonic for the natural invariant 
metric on $X$ coming from the K\"ahler form $\Om_1$.}
$U(N)$-invariant form, given as the trace of 
the Schur functor $\sh_{\l}$ on the curvature
matrix $\Om$. This  was first  demonstrated by Horrocks \cite{H}; 
we   give a different proof in the next section.

\section{Proof of Theorem \ref{mainthm}}
\label{pf}

For each integer $k$ with $1\lequ k\lequ n$, define the {\em special
Schubert class} $\s^k=\s_{(1^k)}$ to be the class
corresponding to the partition $(1^k)=(1,1,\ldots,1)$ of weight $k$.
Over a century ago, Pieri \cite{P} proved the following multiplication 
rule  in $H^*(X)$:
\begin{equation}
\label{pieri}
\s_{\l}\cdot \s^k = \sum \s_{\m},
\end{equation}
summed over all partitions $\m$ whose diagram is obtained by adding
$k$ boxes to the diagram of $\l$, no two in the same row (a more recent
proof is given in \cite{HP}).

A straightforward consequence of   (\ref{pieri}) is the formula
\[
 \s_{\l} = \det (\s^{\l'_i+j-i})_{1\lequ i,j\lequ m},
\]
due to Giambelli \cite{G}, which shows that the special classes
generate the cohomology ring of $G(m,n)$ (the 
corresponding determinantal formula for Schur polynomials was
discovered by Jacobi \cite{J}).
Moreover, one knows that the same
Pieri rule governs the tensor product decomposition
\[
V^{\l}\otimes \wedge^kV = \sum  V^{\m}
\]
of $GL_n$-modules (several different proofs of this are found in
\cite[\S I.6]{FH} and \cite[\S 79]{Z}).  Therefore, Theorem
\ref{mainthm} will follow from the equality $\s^k=
[\Tr(\wedge^k\Om)]$, where $\Tr(\wedge^k\Om)= c_k(Q,h)$ is the {\em
$k$-th Chern form} of $(Q,h)$. Equivalently, we must show that for all
partitions $\l$ and all $k=1,\ldots,n$,
\begin{equation}
\label{integrals}
\int_{X_{\l}} \Tr(\wedge^k\Om)= \delta(\l;\, (1^k)),
\end{equation}
where $\delta(\l;\, \mu)$ is Kronecker's delta.

The integrals (\ref{integrals}) were computed by Chern in \cite[\S 2]{Ch1},
who wrote the integrand using the Maurer-Cartan forms of the unitary
group $U(N)$. The point is that any invariant form on $X=U(N)/(U(m)\times
U(n))$ pulls back to a differential form on $U(N)$, which can be expressed
in terms of the basic invariant forms on $U(N)$.

We will use the differential forms $\om_{ij}$ and $\ov{\om}_{ij}$
%(for $1\lequ i<j\lequ N$)
defined in \cite[\S 5]{T},  of type $(1,0)$ and $(0,1)$, respectively,
which are a scalar multiple of the Maurer-Cartan forms. To describe
them, let $\lie{h}=\{\mathrm{diag}(t_1,\ldots,t_N)\, |\, t_i\in\C\}$
be the Cartan subalgebra of diagonal matrices in $\gl_N$. Consider 
the set of roots 
\[
\Delta=\{t_i-t_j\, |\, 1\lequ i\neq j \lequ N\}\subset \lie{h}^*
\] 
and denote the root $t_i-t_j$ by the pair $ij$. The 
adjoint representation of $\lie{h}$ 
on $\gl_N$ determines a decomposition 
\[
\dis \gl_N=\lie{h}\oplus\sum_{ij\in\Delta}\C\, e_{ij},
\] 
where $e_{ij}$ is the matrix with 1 as the $ij$-th entry and zeroes elsewhere.

Let $\ov{e}_{ij}=-e_{ji}$ and consider the linearly independent
set 
\[ 
B_X=\{e_{ij}, \, \ov{e}_{ij} \  |\  i \lequ m < j \}.
\] 
Extend $B_X$ to a basis $B$ of $\gl_N$ and let $B^*$ be the dual
basis of $\gl_N^*$. For every $e_{ij}$ and $\ov{e}_{ij}$ in $B_X$,
let $\om^{ij}$ and $\ov{\om}^{ij}$ be the corresponding dual basis
vectors in $B^*$, which we shall regard as left
invariant complex one-forms on $U(N)$. Let $\om_{ij}=\gamma\,\om^{ij}$ and
$\ov{\om}_{ij}=\gamma\,\ov{\om}^{ij}$,
where $\gamma$ is a constant such that
$\gamma^2=\frac{i}{2\pi}$.

If $\pi:U(N)\ra X$ denotes the quotient map, then any
smooth form $\eta$ on $X$ pulls back to 
\begin{equation}
\label{forms}
\pi^*\eta=\sum  f_{a_1\ldots a_r b_1\ldots b_s}
\,\omega_{a_1}\wedge\ldots\wedge \omega_{a_r}\wedge
\ov{\omega}_{b_1}\wedge\ldots\wedge\ov{\omega}_{b_s}
\end{equation}
on $U(N)$, with coefficients $f_{a_1\ldots a_r b_1\ldots b_s} \in
C^{\infty}(U(N))$. The differential forms on $U(N)$ which arise in
this way are exactly those which are invariant under the action of the
group $H=U(m)\times U(n)$ ($H$ acts on $C^{\infty}(U(N))$ by right
translation and on $\mathrm{Span}\{\om_{ij}, \, \ov{\om}_{ij} \ |\ i
\lequ m < j \}$ by the dual of the adjoint representation). A smooth
form $\eta$ on $X$ is of $(p,q)$ type precisely when each summand on
the right hand side of (\ref{forms}) involves $p$ unbarred and $q$
barred terms.

The curvature matrix $\Om$ can now be written explicitly (see \cite[\S
4]{GrS} and \cite[Prop.\ 2]{T}); one has
$\Om=\{\Gamma_{\a\b}\}_{m+1\lequ \a,\b \lequ N}$ with
\begin{equation}
\label{eq1}
\Gamma_{\a\b} = \sum_{i=1}^m \om_{i\a}\wedge \ov{\om}_{i\b}.
\end{equation}
Applying the defining expansion of a determinant in terms of its entries, 
we obtain
\begin{equation}
\label{eq2}
\Tr(\wedge^k(\Om))=
\frac{1}{k!}\sum
\mbox{sgn}(\a_1,\ldots,\a_k;\, \b_1,\ldots,\b_k)\, 
\Gamma_{\a_1\b_1}\cdots\Gamma_{\a_k\b_k},
\end{equation}
where the sum is   over all indices $\a_1,\ldots,\a_k,\b_1,\ldots,\b_k$
from the set $\{m+1,\ldots,N\}$, and 
$\mbox{sgn}(\a_1,\ldots,\a_k;\, \b_1,\ldots,\b_k)$ is  zero except when
$(\b_1,\ldots,\b_k)$ is a permutation of $(\a_1,\ldots,\a_k)$, in which
case it equals $+1$ or $-1$ according to the sign of the permutation.
Equations (\ref{eq1}) and (\ref{eq2}) correspond exactly to  
\cite[(12) and (13)]{Ch1}, and   (\ref{integrals}) is proved in
loc.\ cit., Theorem 5, by directly integrating the forms 
(\ref{eq2}) over the Schubert
varieties in $X$ (see also \cite[IV.2]{Ch3}). 

\medskip
\noin
{\bf Remarks.} 1) Although the differential forms on
the right hand side of 
(\ref{eq2}) appear in \cite{Ch1}, their interpretation as 
invariant polynomials in the entries of the curvature matrix  
$\Om=\Omega(Q,h)$
does not. This connection is explained clearly later, in Chern's 
University of Chicago notes \cite[\S 12]{Ch4} and in \cite[\S 8]{Ch5}.

\medskip
\noin
2) For the purposes of the proof, it is not necessary to know that
the ring homomorphisms in (\ref{isom1}) and (\ref{ucw}) are isomorphisms.
Furthermore, each of the three morphisms in (\ref{rho}) maps integral classes
to integral classes. If $\cR_+(GL_n)$ denotes the $\Z$-submodule of 
$\sR_+(GL_n)$ spanned by the polynomial characters
of $GL_n$, we obtain a natural induced homomorphism
\[
\rho_{\Z}: \cR_+(GL_n) \lra H^*(G(m,n),\Z)
\]
which sends $[V^{\l}]$ to the Schubert class $\s_{\l}$ for every partition $\l$.

\section{The problem in other Lie types}
\label{prob}

In this section, we look for an analogue of Theorem \ref{mainthm}
with $GL_n$ replaced by a group of a different Lie type.
Let $G$ be a complex connected reductive Lie group,
$T$ a maximal torus in $G$ and let $\lie{g}$ and $\lie{h}$ be the 
Lie algebras of $G$ and $T$, respectively. 
The adjoint action of $G$ on $\lie{g}$ induces an
action on the ring $\Pol(\lie{g})$ of polynomial functions on $\lie{g}$ with
real coefficients. If $K$ is a maximal compact subgroup of $G$, then
Cartan's theorem \cite{Car} gives a ring isomorphism
\begin{equation}
\label{psiG}
\psi_G:\Pol(\lie{g})^G \ra H^*(BK)
\end{equation}
as before.

If $G$ is not of type A, then extending  (\ref{psiG}) -- on either side -- 
to a sequence analogous to (\ref{rho}) is problematic.
We will go one step further
by using the natural $\l$-ring structure on the representation ring
of $G$. Our main reference for $\l$-rings is \cite{FL}, 
while we learned much of the material that follows from \cite{Fa}, \cite{O},
and \cite{Be}.

A {\em $\l$-ring} is a commutative ring $A$ with a sequence of
operations $\l^i:A\ra A$ such that $\l^0=1$, $\l^1=\mathrm{id}_A$ and
\[
\l^k(x+y)= \sum_{i=0}^k \l^i(x)\l^{k-i}(y)
\]
for all $k\gequ 1$ and for all $x,y \in A$. In addition, we require that 
there are formulas 
\[
\l^k(xy)=P_k(\l^1(x),\ldots,\l^k(x),\l^1(y),\ldots,\l^k(y))
\]
and
\[
\l^k(\l^{\ell}(x))= P_{k,\ell}(\l^1(x),\ldots,\l^{k\ell}(x))
\]
where $P_k$ and $P_{k,\ell}$ are universal polynomials with integer
coefficients, independent of the ring $A$. Setting
$\l_t(x)=\sum_i\l^i(x)\, t^i$, we have $\l_t(x+y)=\l_t(x)\l_t(y)$.
%A homomorphism $f:R_1\ra R_2$ of
%$\l$-rings is a ring homomorphism such that $f(\l^k(x))=\l^k(f(x))$ for all
%$k\gequ 0$.

Let $\cR(G)$ denote the integral representation ring of $G$. Then $\cR(G)$ 
is a $\l$-ring, with 
the $\l$-operations $\l^i:\cR(G)\ra \cR(G)$ induced by the 
exterior powers of representations. There is a unique
$\l$-ring homomorphism $\epsilon:\cR(G)\ra\Z$, called the augmentation,
which associates to each $G$-representation $V$ its dimension $\epsilon(V)$.

Grothendieck \cite{SGA6} introduced $\l$-rings in his work extending
the Hirzebruch-Riemann-Roch theorem to a relative, functorial
setting. He showed that the purely algebraic data of a $\l$-ring $A$
together with an augmentation map $A\ra \Z$ suffice to define
characteristic classes for the elements of $A$.  The Chern classes
constructed in this way take values in the graded ring $\mbox{gr}\,A$
associated to a certain filtration on $A$ coming from the
$\l$-structure, known as the {\em $\gamma$-filtration}. Grothendieck
applied this theory for $A$ equal to the $K$-theory group of vector
bundles on an algebraic variety; we will take $A=\cR(G)$ in the
sequel.

The $\gamma$-operations $\gamma^i$ on $\cR(G)$ are defined by the formula
\[
\gamma_t(x)=\l_{t/(1-t)}(x)=\sum_i\gamma^i(x)\, t^i, \quad \forall x\in \cR(G).
\]
The
$\gamma$-filtration is 
%defined to be 
the decreasing sequence $\{F^k\}_{k\gequ 0}$ where $F^0=\cR(G)$,
$F^1=\Ker(\epsilon)$ and $F^k$ is spanned by the elements
$\gamma^{i_1}(x_1)\cdots\gamma^{i_r}(x_r)$ with $x_1,\ldots,x_r\in
F^1$ and $\sum_{p=1}^r i_p\gequ k$. Let $\mbox{gr}\,
\cR(G)=\bigoplus_{k\gequ 0} F^k/F^{k+1}$ be the associated graded
ring.  For each element $x\in \cR(G)$, there are Chern classes
$c_k(x)$ with values in $\mbox{gr}^k \cR(G)$.  By definition,
$c_k(x)=\gamma^k(x-\epsilon(x))$.

\medskip
\noin
{\bf Example 4.} Suppose that the torus $T$ has rank $n$, and let  $\Lambda$
be its character group. We can identify $\Lambda$ with the multiplicative
group of monomials $\a_1^{m_1}\cdots\a_n^{m_n}$ with $m_i\in \Z$, where
$\a_i^{\pm 1}$ corresponds to the map $(t_1,\ldots,t_n)\mapsto t_i^{\pm 1}$,
for $1\lequ i \lequ n$. We then have
\[
\cR(T)=\Z[\Lambda]=\Z[\a_1,\a_1^{-1},\ldots,\a_n,\a_n^{-1}].
\]
The relations $\l_t(\xi)=1+\xi t$ and $\epsilon(\xi)=1$ for $\xi\in
\Lambda$ determine the $\l$-structure and augmentation on $\cR(T)$. It
is straightforward to check that the $\gamma$-filtration of $\cR(T)$
coincides with its $F^1$-adic filtration, that is, $F^k=(F^1)^k$ for
all $k\gequ 1$. The map $\xi\mapsto \xi -1$ gives a canonical
isomorphism from $\Lambda$ to the additive group $\mbox{gr}^1 \cR(T)=
F^1/F^2$. If $u_i\in \mbox{gr}^1 \cR(T)$ denotes the image of $\a_i$
under this map, then the elements $u_1,\ldots,u_n$ are algebraically
independent over $\Z$, and we have $\mbox{gr}\, \cR(T)=
\Z[u_1,\ldots,u_n]\cong \mbox{Sym}(\Lambda)$ (see e.g.\ \cite[Prop.\
3.2]{O}). We may thus identify the $k$-th Chern class of $\xi_1+\ldots
+\xi_r$ (where $\xi_i\in\Lambda$) with the element
$e_k(\xi_1,\ldots,\xi_r) \in \mbox{Sym}(\Lambda)$.

\medskip

The Weyl group $W$ acts on $T$, hence on the ring $\cR(T)$, and the
restriction map $\eta:\cR(G)\ra \cR(T)$ induces an isomorphism of
$\lambda$-rings $\cR(G)\ra \cR(T)^W$.  Let us now pass to real
coefficients: define $\sR(G)=\cR(G)\otimes\R$, $\mbox{gr}\,
\sR(G)=\mbox{gr} \,\cR(G)\otimes \R$ and similarly for the torus $T$.
One can then show that the map $\eta$ respects the
$\gamma$-filtrations of both $\sR(G)$ and $\sR(T)$, and thus induces a
graded ring homomorphism $\mbox{gr}(\eta) : \mbox{gr}\, \sR(G) \ra
\mbox{gr}\, \sR(T)$, which maps $\mbox{gr}\, \sR(G)$ isomorphically
onto the invariant subring $(\mbox{gr}\, \sR(T))^W$ (see \cite{Fa},
\cite{O}, and \cite{Be} for further discussion and proofs).

Assume that the Lie groups $G$, $T$ and their Lie algebras are defined over the
real numbers, and view $\lie{h}$ as a vector space over $\R$.
Applying the
map which takes a character in $\Lambda$ to its 
derivative in $\lie{h}^*$, we can identify $(\mbox{gr}\, \sR(T))^W
\cong (\mbox{Sym}(\Lambda))^W\otimes \R$ with the algebra
$\Pol(\lie{h})^W$ of $W$-invariant polynomial functions on $\lie{h}$
(again with real coefficients). Chevalley proved
that the restriction homomorphism $\Pol(\lie{g})\ra\Pol(\lie{h})$ maps 
$\Pol(\lie{g})^G$ isomorphically onto $\Pol(\lie{h})^W$. Combining this 
with the preceeding ingredients shows that
(\ref{psiG}) extends to a sequence of isomorphisms
\[
\mbox{gr}\,\sR(G)\lra\Pol(\lie{g})^G \stackrel{\psi_G}\lra H^*(BK)
\]
which respects the natural grading in all three rings.

\medskip
\noin
{\bf Example 5.}
Let $x$ be the element of $\sR(GL_n)$ which corresponds to the standard
representation of $GL_n(\C)$ on $\C^n$. It is easy to see that the $k$-th
Chern class $c_k(x)$ (in the above sense) is the invariant polynomial 
$A\mapsto \Tr(\wedge^kA)$ on $\lie{g}\lie{l}_n$. Theorem \ref{mainthm}  
shows that the Chern classes of the standard representation of $GL_n$
map naturally to the special Schubert classes in $G(m,n)$. 
The fact that the target ring $\mbox{gr}\,\sR(GL_n)$ for the $\sR(GL_n)$
Chern classes is naturally graded 
isomorphic to $\sR_+(GL_n)$ is a type A phenomenon.

\medskip
\noin
{\bf Example 6.} Let $G=Sp_{2n}(\C)$ be the symplectic group of rank $n$ and
$K=Sp(2n)$ the homonymous compact subgroup. Let $y\in \sR(Sp_{2n})$ 
be the class of the standard representation of $G$ on $\C^{2n}$.
In this case the $2k$-th Chern class $c_{2k}(y)$ is the invariant polynomial
$A\mapsto \Tr(\wedge^{2k}A)$ on $\lie{g}=\lie{s}\lie{p}_{2n}$, 
and these classes for $1\lequ k \lequ n$
generate the algebra $\Pol(\lie{g})^G$ (the odd Chern classes of
$y$ vanish). Observe that $\Pol(\lie{g})^G$
 is isomorphic to $\Pol(\lie{g}\lie{l}_n)^{GL_n}$, up to a 
doubling of degrees.

The classifying space $BSp(2n)$ may be identified with the infinite
quaternionic Grassmannian $G_{{\mathbb H}}(\infty,n)$, which is the
inductive limit of the finite Grassmannians $Y=G_{{\mathbb H}}(m,n)$
over all $m>0$. Note that $Y=Sp(2N)/(Sp(2m)\times Sp(2n))$ is a
compact oriented manifold of real dimension $4mn$, which parametrizes
$m$-dimensional left ${\mathbb H}$-linear subspaces of ${\mathbb
H}^N$. There is a rank $n$ quotient bundle $Q\ra Y$ which is a
pullback of the universal quotient bundle over $G_{{\mathbb
H}}(\infty,n)$. The cell decomposition and cohomology ring of $Y$ are
identical to those of the complex Grassmannian $G(m,n)$, again up to a
doubling of degrees (see \cite{Bo1}, \cite[Vol.\ III, Chp.\ XI]{GHV},
and \cite[Appendix A]{PR}).  In particular, for each partition $\l$ as
in Sec.\ \ref{svcc}, there is a `Schubert variety' $Y_{\l}$ in $Y$,
defined by the same inequalities that define $X_{\l}$ in $G(m,n)$, and
a `Schubert class' $\t_{\l}=[Y_{\l}]\in H^{4|\l|}(Y)$. By arguing as
in Sec.\ \ref{pf}, one can show that the induced ring homomorphism
\[
\mbox{gr}\, \sR(Sp_{2n}) \lra \Pol(\lie{s}\lie{p}_{2n})^{Sp_{2n}} 
\stackrel{\psi_{Sp}}\lra H^*(BSp(2n))\lra H^*(G_{{\mathbb H}}(m,n))
\]
maps the generators $c_{2k}(y)$ to the `special Schubert classes'
$\t_{(1^k)}$. 

We remark that the irreducible characters $\chi^{\l}$ of
$Sp_{2n}$  are parametrized by partitions of length at most $n$, as is
the case for the polynomial characters of $GL_n$. 
Moreover, in the product decomposition
\[
\chi^{\l}\chi^{\mu}=\sum_{\n}c^{\n}_{\l\m}\,\chi^{\nu},
\]
the structure constants $c^{\n}_{\l\m}$ agree with those in (\ref{gltensor}),
whenever $|\n|=|\l|+|\m|$ (see e.g.\ \cite[Prop.\ 2.5.2]{KT}). 

\medskip

In the case of the orthogonal groups $G=O_n(\C)$, the classifying space
$BO(n)$ is an infinite real Grassmannian, which is a limit of finite
dimensional real Grassmannians $G_{\R}(m,n)$. However there is no
clear analogue of $GL_n(\C)$ Schubert calculus on these manifolds, some
of which are not even orientable.

\section{Concluding remarks}

In this final section, we briefly discuss some of the early works
where the representation theory of Lie groups was applied to
study the cohomology of  homogeneous spaces, and which 
relate somehow to the present paper. For more details
about the early investigations of the Chern-Weil circle of ideas and its  
applications, we recommend Chern's address at the 1950
International Congress \cite{Ch2}, Weil's letters, written in 1949 
and first published in \cite{W}, 
the survey articles by Samelson \cite{Sa} 
and Borel \cite{Bo2}, and the historical notes in \cite[Vol.\ III]{GHV}.

In the middle of the last century, there emerged two effectively different 
approaches to the study of the cohomology of principal bundles and
homogeneous spaces of Lie groups. 
The techniques used in this paper derive from \'E.\ Cartan's method of
invariant differential forms \cite{ECar} and its later extension by H.\ Cartan 
\cite{Car}. An alternative approach, espoused by Borel \cite{Bo1} among
others, used the methods of classical algebraic topology. The theorems
proved by these two schools were frequently identical, although the results
obtained by topological methods were usually with $\Z$ or $\Z_p$ coefficients.
These latter techniques were applied by Borel and Hirzebruch \cite{BH} to
connect the theory of characteristic classes to the cohomology of 
homogeneous spaces $G/U$, by interpreting the former as elementary
symmetric functions in certain roots of $G$ (or their squares).

Ehresmann \cite{E} investigated the topology of complex Grassmann
manifolds (and other hermitian symmetric spaces) by studying the
algebra of $K$-invariant differential forms on them ($K=U(N)$ for
$X=G(m,n)$).  This relies on the fact that the invariant forms are
harmonic for the natural hermitian structure on $X$, which implies
that the ring of all such forms is isomorphic to $H^*(X)$. Kostant
\cite{Ko1} \cite{Ko2} later found analogues of these results for
arbitrary (generalized) flag manifolds.  The representation theory
used to determine the $K$-invariant forms in this program does not
directly relate the multiplicities $c^{\nu}_{\l\mu}$ in equations
(\ref{gltensor}) and (\ref{gmncup}). Note however that the cited works
of \'E.\ Cartan and Ehresmann were used by Chern in his fundamental
paper on the characteristic classes of complex manifolds
\cite{Ch1}. More recently, Stoll \cite{St} used fiber integration 
to study the algebra of invariant forms on the Grassmannian, 
but his work does not address the question posed in the
Introduction.

%For $G$ as in Section \ref{prob}, there is a Chern character homomorphism
%\[
%ch: \sR(G)\lra \prod_{i\gequ 0} \mbox{gr}^i\sR(G)
%\]
%which extends to a ring isomorphism 
%\[
%\wh{ch} : \wh{\sR(G)} \stackrel{\cong}\lra \prod_{i\gequ 0} \mbox{gr}^i\sR(G),
%\]
%where $\wh{\sR(G)}$ denotes the completion of $\sR(G)$ with respect
%to the $\gamma$-filtration (see \cite[Lemma 3.2]{O} and also 
%\cite[\S 4.8]{AH} 
%for a topological interpretation of $\wh{\sR(G)}$). The induced maps
%from $\sR(G)$ and $\wh{\sR(G)}$ to $\wh{\Pol}(\lie{g})^G$ are not 
%suited to our purposes.

Following \cite{SGA6}, \cite{V} and \cite{Be}, 
the isomorphism between $\mbox{gr}\sR(G)$\ and $\Pol(\lie{g})^G$
in Sec.\ \ref{prob} may be used to construct the Chern-Weil (or characteristic) 
homomorphism
in algebraic geometry. Let $P\ra X$ be a principal $G$-bundle over a smooth
algebraic variety $X$ and let $CH^*(X)$ denote the Chow group of algebraic cycles
on $X$ modulo rational equivalence.
The Grothendieck group $K(X)$ of vector bundles on $X$ is a 
$\l$-ring, with the $\l$-operations induced by exterior powers.
According to \cite[Exp.\ XIV]{SGA6}, the graded ring
$\mbox{gr}K(X)\otimes\R$ is canonically isomorphic to the real Chow
ring $CH^*_{\R}(X)=CH^*(X)\otimes\R$.
There is a natural $\l$-ring homomorphism
$\pi:\cR(G)\ra K(X)$, defined by sending a representation $G\ra GL(E)$ to 
the associated vector bundle $P\times_{G}E$ over $X$.  The characteristic homomorphism is the induced map
\[
\mbox{gr}(\pi)_{\R} : \Pol(\lie{g})^G\lra CH_{\R}^*(X).
\]

\end{document}